\documentclass[12pt,twoside]{article} 

\usepackage[usenames]{color}
\usepackage[colorlinks=false,
linkcolor=webgreen,
filecolor=webbrown,
citecolor=webgreen]{hyperref}

\definecolor{webgreen}{rgb}{0,.5,0}
\definecolor{webbrown}{rgb}{.6,0,0}

\date{} 

\newcommand{\li}{\mathop{\mathrm{li}}}
\newcommand{\seqnum}[1]{\href{http://oeis.org/#1}{\textcolor{black}{\underline{#1}}}}

\begin{document} 

\centerline{\bf International Mathematical Forum, Vol.\,10, 2015, no.\,6, 283--288} 


\centerline{\bf http://dx.doi.org/10.12988/imf.2015.5322}

\centerline{} 

\centerline{} 

\centerline {\Large{\bf Verification of the Firoozbakht Conjecture}} 

\centerline{} 

\centerline{\Large{\bf for Primes up to Four Quintillion}} 

\centerline{} 

\centerline{\bf {Alexei Kourbatov}} 

\centerline{} 

\centerline{www.JavaScripter.net/math} 

\centerline{15127 NE 24th Street \#578} 

\centerline{Redmond, WA, USA} 


\centerline{} 







\newtheorem{Theorem}{\quad Theorem}[section] 

\newtheorem{Definition}[Theorem]{\quad Definition} 

\newtheorem{Corollary}[Theorem]{\quad Corollary} 

\newtheorem{Lemma}[Theorem]{\quad Lemma} 

\newtheorem{Example}[Theorem]{\quad Example} 

\centerline{}

{\footnotesize Copyright $\copyright$ 2015 Alexei Kourbatov. This is an open access article distributed under the Creative Commons Attribution License, which permits unrestricted use, distribution, and reproduction in any medium, provided the original work is properly cited.}

\begin{abstract} 
If $p_k$ is the $k$th prime, the Firoozbakht conjecture states that the sequence $(p_k)^{1/k}$ 
is strictly decreasing. We use the table of first-occurrence prime gaps in combination with 
known bounds for the prime-counting function to verify the Firoozbakht conjecture 
for primes up to four quintillion ($4\times10^{18}$).
\end{abstract} 



\section{Introduction}
We will examine a conjecture that was first stated in 1982 by 
the Iranian mathematician Farideh Firoozbakht from the University of Isfahan \cite{rivera}. 
It appeared in print in {\it The Little Book of Bigger Primes} by Paulo Ribenboim \cite[p.\,185]{ribenboim}. 
The statement is as follows:

\medskip\noindent
{\bf Firoozbakht's Conjecture.}
If $p_k$ is the $k$th prime, then the sequence $(p_k)^{1/k}$ is strictly decreasing.
Equivalently, for all $k\ge1$ we have
\begin{equation}\label{ineqf}
p_{k+1}^k < p_k^{k+1}.
\end{equation}
The Firoozbakht conjecture is one of the strongest upper bounds for prime gaps.
As we will see from Table 1 below, it is somewhat stronger than {\it Cram\'er's conjecture} 
proposed about half a century earlier by the Swedish mathematician Harald Cram\'er \cite{cram}:

\medskip\noindent
{\bf Cram\'er's Conjecture.}  If $p_{k}$ and $p_{k+1}$ are consecutive primes, then we have
$p_{k+1}-p_{k} = O(\log^2 p_{k})$ or, more specifically, 
$$
\limsup_{k \to \infty} {p_{k+1}-p_{k} \over \log^2p_k} ~=~ 1.
$$
For the sake of numerical comparison with (\ref{ineqf}), let us use a modified form of 
Cram\'er's conjecture stated below. 

\medskip\noindent
{\bf Modified Cram\'er Conjecture.} If $p_{k}$ and $p_{k+1}$ are consecutive primes, then 
\begin{equation}\label{ineqc}
p_{k+1}-p_{k} ~<~ \log^2 p_{k+1}.
\end{equation}
This modified form allows us to make predictions of an upper bound 
for any given prime gap; Table 1 lists a few examples of such upper bounds.
{\footnotesize
\begin{center}TABLE 1 \\
 Prime gap bounds predicted by the modified Cram\'er and Firoozbakht conjectures \\[0.5em]
\begin{tabular}{rrrrr}
\hline
   &\multicolumn{2}{c}{  }
   &\multicolumn{2}{c}{Upper bounds for $p_{k+1}$, as predicted by: \phantom{\fbox{$1^1$}} } \\
$k$&\multicolumn{2}{r}{Consecutive primes } &{Modified Cram\'er conjecture} & {Firoozbakht conjecture} \\
   & $p_k$  & $p_{k+1}$ &(solution of $x=p_k+\log^2x$) & (solution of $x^k=p_k^{k+1}$) \\
[0.5ex]\hline
\vphantom{\fbox{$1^1$}}
      5 &      11 &      13 &      19.964 &      17.769 \\
     26 &     101 &     103 &     124.255 &     120.618 \\
    169 &    1009 &    1013 &    1057.493 &    1051.152 \\
   1230 &   10007 &   10009 &   10091.999 &   10082.220 \\
   9593 &  100003 &  100019 &  100135.579 &  100123.090 \\
  78499 & 1000003 & 1000033 & 1000193.874 & 1000179.012 \\
 664580 &10000019 &10000079 &10000278.794 &10000261.534 \\
\hline
\end{tabular}
\end{center}
}
Table 1 shows that, given $k$ and $p_k$,  
the Firoozbakht conjecture (\ref{ineqf}) yields a tighter bound for $p_{k+1}$ than 
the modified Cram\'er conjecture (\ref{ineqc}). Indeed, the Firoozbakht upper bound (last column)
is below the Cram\'er upper bound by approximately $\log p_k$. 
In {\it Cram\'er's probabilistic model of primes} \cite{cram,kourbatov2014}
the parameters of the distribution of maximal prime gaps
suggest that inequalities (\ref{ineqf}) and (\ref{ineqc}) are both true with probability 1;
that is, almost all\footnote{
In Cram\'er's model with $n$ urns, the limiting distribution of maximal ``prime gaps''
is the Gumbel extreme value distribution with scale $a_n\sim n/\li n = O(\log n)$
and mode $\mu_n = n\log(\li n)/\li n +O(\log n) = \log^2 n - \log n\log\log n +O(\log n)$ 
\cite[OEIS \seqnum{A235402}]{kourbatov2014,oeis}; here $\li n$ denotes the logarithmic integral of $n$.
Hence, for large $n$, all maximal gap sizes are below $\log n(\log n-1)$, 
except for a vanishing proportion of maximal gaps.
}  
maximal prime gaps in Cram\'er's model satisfy (\ref{ineqf}) and (\ref{ineqc}).
One may take this as an indication that any violations of (\ref{ineqf}) and (\ref{ineqc})
occur exceedingly rarely (if at all).

\section{Computational verification for small primes}

When primes $p_k$ are not too large, one can directly verify inequality (\ref{ineqf}) by computation.
A simple program that takes a few seconds to perform the verification for $p_k<10^6$ is available on the author's website.
The program outputs the numeric values of $k$, $p_k$, $p_k^{1/k}$, and an {\tt OK} if the value
of $p_k^{1/k}$ decreases from one prime to the next.
It will output {\tt FAILURE} if the conjectured decrease 
does not occur. Result: all {\tt OK}s, no {\tt FAILURE}s. Here is a sample of the output:
{\small
\begin{verbatim}
      k       p      p^(1/k)      OK/fail 
      1       2     2.000000000    OK    
      2       3     1.732050808    OK    
      3       5     1.709975947    OK    
      4       7     1.626576562    OK    
      5      11     1.615394266    OK     
      6      13     1.533406237    OK     
      7      17     1.498919872    OK    
      8      19     1.444921323    OK     
      9      23     1.416782203    OK     
     10      29     1.400360331    OK     
     11      31     1.366401518    OK    
     12      37     1.351087503    OK 
    ...     ...     ...
  78494  999953     1.000176022    OK   
  78495  999959     1.000176020    OK   
  78496  999961     1.000176018    OK   
  78497  999979     1.000176016    OK   
  78498  999983     1.000176014    OK   
  78499 1000003     1.000176012    OK    
  78500 1000033     1.000176010    OK   
\end{verbatim}
}

\noindent 
Thus primes $p_k<10^6$ do not violate (\ref{ineqf}). What about larger primes?

\section{What if we do not know $k=\pi(p_k)$?}

For large primes $p_k$, the exact values of $k=\pi(p_k)$ are not readily available.
Nevertheless, the Firoozbakht conjecture can often be verified in such cases too.
When we do not know $\pi(p_k)$ exactly, we can use these bounds for the
prime-counting function $\pi(x)$:
\begin{equation}\label{pi11}
\pi(x) < {x\over\log x - 1.1} ~\mbox{ for } x\ge 60184 \ \qquad
\mbox{ \cite[p.\,9, Theorem 6.9]{dusart} } \qquad\qquad 
\end{equation}
\begin{equation}\label{pi12}
\pi(x) < {x\over\log x - 1.2} ~\mbox{ for } x\ge 4 
\quad\mbox{\small (from (\ref{pi11}) + computer check for $x<10^5$). }
\end{equation}
Taking the log of both sides of (\ref{ineqf}) and rearranging, we find that the Firoozbakht conjecture 
(\ref{ineqf}) is equivalent to
\begin{equation}\label{firlog}
\pi(p_k) < {\log p_k \over \log p_{k+1} - \log p_k}.
\end{equation}
If we know $p_k$ and $p_{k+1}$ (where $p_k>60184$) but do not know $\pi(p_k)$, 
then instead of (\ref{firlog}) we may check the stronger condition
\begin{equation}\label{combineq}
{p_k\over\log p_k - 1.1} < {\log p_k \over \log p_{k+1} - \log p_k}.
\end{equation}
For a larger range of applicability ($p_k>4$)\footnote{
regardless of the computation of Sect.\,2 which already proves (\ref{ineqf}) for all $p_k<10^6\ldots$
}
we may check another (still stronger) condition:
\begin{equation}\label{combineq2}
{p_k\over\log p_k - 1.2} < {\log p_k \over \log p_{k+1} - \log p_k}.
\end{equation}
If (\ref{combineq}) or (\ref{combineq2}) is true, so is (\ref{firlog});
and the exact value of $\pi(p_k)$ is not needed to check (\ref{combineq}), (\ref{combineq2}).

\section{Verification for all gaps of a given size $g$}\label{sect4}

We will take (\ref{combineq}) and (\ref{combineq2}) one step further and make $p_k$ a variable ($x$); 
then $p_{k+1}=x+g$, where $g$ is the gap size. We can now solve the resulting simultaneous inequalities
\begin{equation}\label{combineqx}
0 < {x\over\log x - 1.1} < {\log x \over \log(x+g) - \log x} \quad\mbox{ with }\quad x>60184,
\end{equation}
or, if we are interested in a larger range of applicability,
\begin{equation}\label{combineqx2}
0 < {x\over\log x - 1.2} < {\log x \over \log(x+g) - \log x} \quad\mbox{ with }\quad x>4. \phantom{1111}
\end{equation} 
Here we use the gap size $g$ as a parameter. 
In combination with a table of first-occurrence prime gaps \cite{nicely}, 
the solution of (\ref{combineqx}) and/or (\ref{combineqx2}) will tell us whether a prime gap 
of size $g$ may violate the Firoozbakht conjecture for primes $p_k\approx x$. 
Consider the following examples.

\medskip\noindent
{\bf Example 1.} {\it Can a prime gap of size $150$ violate the Firoozbakht conjecture?} \\
To answer this question, we substitute $g=150$ into (\ref{combineqx}),
$$
0 < {x\over\log x - 1.1} < {\log x \over \log(x+150) - \log x} \quad\mbox{ with }\quad x>60184,
$$
solve for $x$ and find the ``safe bound''
$$
x\ge365323 \quad(\mbox{or, more precisely, } x>365322.7038);
$$
that is, a gap of 150 does not violate (\ref{ineqf}) if such a gap occurs between primes above 365323.
But there are no prime gaps of size 150 below 365323; in fact, the table of 
first-occurrence prime gaps \cite{nicely} indicates that the first such gap follows the prime 13626257. 
Therefore, a prime gap of size 150 can never violate the Firoozbakht conjecture (\ref{ineqf}).

\medskip\noindent
{\bf Example 2.} {\it Can a prime gap of size $2$ (twin primes) violate the Firoozbakht conjecture?} 
\ \ We substitute $g=2$ into (\ref{combineqx2}), 
$$
0 < {x\over\log x - 1.2} < {\log x \over \log(x+2) - \log x} \quad\mbox{ with }\quad x>4,
$$
solve for $x$ and find the ``safe bound''
$$
x\ge8 \quad(\mbox{or, more precisely, } x>7.8745).
$$
So a prime gap of size 2 does not violate (\ref{ineqf}) if this gap occurs between primes above 8.
But we already know that gaps between primes below 8 do not violate the conjecture either (Section 2). 
Therefore, a prime gap of size 2 (twin primes) can never violate (\ref{ineqf}).

\medskip\noindent
We have repeated the computation of the above examples 
for all even values of the gap~size $g\in[2,1476]$ and found that 
{\it none of these gap sizes could possibly violate} (\ref{ineqf}).
A tabulation of ``safe bounds'' by gap size is available on the author's website \cite{safebounds}.
For $g=2$ and $4$, we had to manually check (\ref{ineqf}) for a couple of gaps of size $g$
between primes below the respective safe bounds.
For $g\in[6,1476]$, the actual first occurrence of prime gap $g$ is already safe. 
(``Close calls'' occur for record prime gaps in OEIS sequence \seqnum{A005250} \cite{oeis}.)
From the prime gaps table \cite{nicely} we also know that gaps larger than 1476 do not occur 
below $4\times10^{18}$. Thus the validity of Firoozbakht conjecture (\ref{ineqf}) 
has been verified for primes up to $4\times10^{18}$. We have obtained the following theorem:

\medskip\noindent
{\bf Theorem } {\it Inequality (\ref{ineqf}) is true for all primes} $p_k < 4\times10^{18}$.
%

\bigskip\noindent
{\bf Acknowledgements.} 
The author expresses his gratitude to all contributors and editors of the websites 
{\it OEIS.org} and {\it PrimePuzzles.net}, especially to Farideh Firoozbakht 
for proposing a very interesting conjecture.
Thanks are also due to Pierre Dusart for proving the $\pi(x)$ bound (\ref{pi11}),
to Tom\'as Oliveira e Silva, Siegfried Herzog, and Silvio Pardi 
whose computation extended the table of first-occurrence prime gaps \cite{toes2014},
and to Thomas Nicely for maintaining the said table in an easily accessible format \cite{nicely}.

\pagebreak
\footnotesize{
\noindent
{\bf Endnotes.}
All first-occurrence gaps between primes $p_k<2^{64}$ have been found by 2019. \linebreak
The results are reported at {\it Mersenneforum.org}\, and included in T.\,R.~Nicely's tables
\cite{nicely}, as well as in OEIS \cite[\seqnum{A014320}, \seqnum{A335366}, \seqnum{A335367}]{oeis}.
Using the method of Section \ref{sect4}, these results allow us to confirm that
Firoozbakht's conjecture (\ref{ineqf}) holds for all primes $p_k<2^{64}$.
As in Example 1, the WolframAlpha command
$$
\mbox{\tt solve 0 < x/(ln(x)-1.1) < ln(x)/(ln(x+1918)-ln(x))}
$$
uses gap $g=1918$ and yields a safe bound a little below $2^{64}$. Even lower safe bounds are
obtained with the same command for gaps $g<1918$ --- and the actual first occurrences of gaps $g$
are already safe (i.e., the corresponding primes $p_k$ are larger than their respective safe bounds).
Thus, extended computations have shown that 
\begin{itemize}
\item Firoozbakht's conjecture (\ref{ineqf}) is true for all primes $p_k<2^{64}$;
\item prime gaps of size $g<1920$ cannot violate (\ref{ineqf}).
\end{itemize}
{\it (Endnotes added 5 Jan 2023. The URLs \cite{nicely} updated.)}

}



\begin{thebibliography}{99} 

\bibitem{cram}
H.~Cram\'er, On the order of magnitude of the difference between consecutive prime numbers. 
{\it Acta Arith.} {\bf 2} (1936), 23-46. 

\bibitem{dusart}
P.~Dusart, Estimates of some functions over primes without R.H., 
arXiv:1002.0442 (2010). \url{http://arxiv.org/abs/1002.0442}

\bibitem{kourbatov2014}
A.~Kourbatov, The distribution of maximal prime gaps in Cram\'er's probabilistic model of primes.
{\it Int.~Journal of Statistics and Probability} {\bf 3} (2014), No.\,2, 18-29.
arXiv:1401.6959. \url{http://arxiv.org/abs/1401.6959}

\bibitem{safebounds}
A.~Kourbatov, Checking the Firoozbakht conjecture: safe bounds, 2015. \\ 
\url{http://www.javascripter.net/math/primes/firoozbakhtconjecturebounds.htm}

\bibitem{nicely}
T.~R.~Nicely, First occurrence prime gaps, preprint, 2014 and later years. \\
\url{https://faculty.lynchburg.edu/~nicely/gaps/gaplist.html} \\
\url{http://oeis.org/A000101/a000101.pdf}

\bibitem{toes2014}
T.~Oliveira e Silva, S.~Herzog, and S.~Pardi,
Empirical verification of the even Goldbach conjecture and computation of prime gaps up to $4\cdot10^{18}$, 
{\it Math.~Comp.} {\bf 83} (2014), 2033-2060. \ 
\url{http://dx.doi.org/10.1090/s0025-5718-2013-02787-1}

\bibitem{ribenboim}
P.~Ribenboim, {\it The Little Book of Bigger Primes}, %
Springer, New York, NY, 2004. 
\url{http://dx.doi.org/10.1007/b97621 }

\bibitem{rivera} 
C. Rivera (ed.), Conjecture 30. The Firoozbakht Conjecture, 2002. \\
\url{http://www.primepuzzles.net/conjectures/ } 

\bibitem{oeis}
N.~J.~A.~Sloane (ed.), {\it The On-Line Encyclopedia of Integer Sequences}, 2014.
Published electronically at \url{http://oeis.org/}. Sequences
\seqnum{A014320},
\seqnum{A235402},
\seqnum{A335366},
\seqnum{A335367}.

\end{thebibliography}
\end{document}